# ENDS OF ITERATED FUNCTION SYSTEMS

GREGORY R. CONNER AND WOLFRAM HOJKA

ABSTRACT. Guided by classical concepts, we define the notion of *ends* of an iterated function system and prove that the number of ends is an upper bound for the number of nondegenerate components of its attractor. The remaining isolated points are then linked to idempotent maps. A commutative diagram illustrates the natural relationships between the infinite walks in a semigroup and components of an attractor in more detail. We show in particular that, if an iterated function system is one-ended, the associated attractor is connected, and ask whether every connected attractor (fractal) conversely admits a one-ended system.

## 1. INTRODUCTION

The main aim of the current article is to offer a conceptually new treatment to the study of a general iterated function system (IFS) in which one has no *a priori* information concerning the attractor. We describe how purely algebraic properties of iterated function systems can anticipate the topological structure of the corresponding attractor. Our method convolves four natural viewpoints: algebraic, geometric, asymptotic, and topological. Algebraically, the semigroup of an IFS, which consists of all countably many compositions of the functions in the system (most notably the idempotent elements), carries a surprising amount of information about the attractor. The arguments are basically geometric in nature. Infinite walks, which reside in the intersection of algebra and geometry, play an important role. The asymptotic point of view is embodied by the classical concept of *ends*. Topologically, we relate components of the attractor and especially the isolated points (the degenerate components) to the algebraic and asymptotic properties of the system.

Previous literature has described the attractor of an IFS by detailed analysis of the involved maps. For example [Hat85], [KL00], and more recently, [AT04], or [LT06] require a minute understanding of the action of the IFS on the attractor. In consonance with related fields including geometric group theory and dynamical systems our ambition is to investigate the attractor by examining the relations in the system only on a formal, algebraic level.

Accordingly, we obtain an asymptotic decomposition of the Cayley graph in Theorem 8. This decomposition yields a bound on the number of components of the attractor by the number of ends of the IFS, and identifies isolated points as idempotents in the semigroup, see Theorems 3 and 10. Of special interest is the case where there exist no idempotents and only one end, as that implies the attractor to be connected. Theorem 13 offers another way to characterize this situation by the connectivity of a different, finite graph. This has some superficial similarity to Hata's criterion in [Hat85], but whereas his conditions are topological, ours concern algebraic relations. Conversely we demonstrate how a number of well-known classical connected fractals can be represented by one-ended systems. Encouraged by this, Conjecture 4 asks if *every* connected attractor can be realized by a one-ended IFS.

Diagram 11 testifies to the naturality of the presented connections between the algebraic and the topological structure of the semigroup and attractor respectively. It highlights a dichotomy within the set of infinite walks in the semigroup which is reflected in the topology of the attractor. On one hand, the coarse geometry of the walks is reflected in the language of ends and correspondingly in the nondegenerate components of the attractor. On the other hand, the local, finite structure of cycling walks (giving rise to dead-ends in the graph) manifests itself in idempotents in the semigroup (Lemma 9).

While isolated points in a given attractor guarantee, by virtue of Theorem 10, the existence of idempotents in any corresponding semigroup, Example 19 shows that the converse is not true. Furthermore, one cannot preclude the existence of idempotents in the semigroup by looking for idempotents only among the functions which generate the IFS, as the Examples 20 and 21 illustrate.

The literature contains other interesting methods to describe an attractor such as in [Rud92] (by means of inverse limits) and for the special case of the Sierpiński triangle in [Kai03], again employing a significantly more topological approach. Another article that deals with that space is [ADTW09]. In some sense their technique is the opposite of ours: knowing that the Sierpiński triangle has a dense set of local cutpoints the geometry of the cutpoint set is used to describe paths in the space via a coordinate system based on the word structure of the IFS.

## 2. Main Results

There is a standard construction that encodes the algebraic structure of a group into a graph, first described by Arthur Cayley in [Cay78], which has turned out to be one of the most important tools in low dimensional topology and geometric group theory. The *Cayley graph* is the geometric model of a discrete group, using lengths of words as a measure of distance. We describe the natural analogue of the Cayley graph for iterated functions systems. While the geometric realization of a group is its Cayley graph, the Cayley graph of an IFS only hints at the geometry of the attractor.

Consider an iterated function system $F$ consisting of finitely many contracting maps from a complete metric space $M$ to itself. Its attractor (*fractal*, see [Hut81]) $A$ is the uniquely determined compact subspace of $M$ with $\bigcup FA = \bigcup \{f(A) : f \in F\} = A$. Composition of functions yields a semigroup structure on $M^M$, the set of maps from $M$ to itself.

**Definition 1.** Denote with $S := S(F)$ the subsemigroup generated by the finitely many maps in $F$. Associated to this semigroup is its *Cayley graph* $\Gamma(S)$, with vertices $S$ and edges $(s, f) \in S \times F$ connecting $s$ and $s \circ f$.

This graph is the central object of our study. Note that its structure depends only on the *equality* of compositions of functions in the system, not on any analytical, geometric, or topological properties they may or may not possess. Thus this is a purely algebraic object and so has only passing similarity to other constructs in the literature such as that of [Hat85].

Based on work of Freudenthal [Fre31], the set of ends can be defined for a generic topological space $X$ to be the inverse limit of the components of the complement of a compact subset of $X$, where the bonding maps are induced by inclusion. In the more restricted setting of graphs the following characterization is more convenient (cf. also [Hal64]):

**Definition 2.** A sequence of vertices in a graph $\Gamma$ is called a *walk* if two consecutive vertices are connected by an (oriented) edge in the graph. If additionally all vertices are distinct, it is called a *ray*. Two rays are *equivalent* if for every finite set $E$ of vertices in $\Gamma$ they have subsequences that lie in a common path component of $\Gamma \smallsetminus E$. These equivalence classes form the *ends* of the graph.

An IFS is said to have the same set of ends as its Cayley graph.

**Remark.** Whilst not the case for generic graphs, orientation has no bearing on the number of ends of an IFS. In particular, a graph has *one end* if the complement of any finite set of vertices has exactly one infinite connected component. A graph has *zero ends* if it is finite.

We ask the reader for forgiveness in that we deviate from the standard notation in the most degenerate of all cases. We will consider the single point of a one point space as nonisolated. This allows us to state the main theorem:

**Theorem 3.** *The number of ends of an iterated function system is an upper bound for the number of nondegenerate components in the attractor. In addition, there may be countably many isolated points, each corresponding to an idempotent in the semigroup generated by the system.*

*Consequently, the attractor is connected if the system has one end and the generated semigroup has no idempotent elements.*

As the examples in the next section show, an IFS whose attractor is connected need not be one-ended. On the other hand, it is natural to ask if a system can be *augmented* so that the above statement can be reversed:

**Conjecture 4.** *If the attractor $A$ of an IFS $F$ is connected, then there is a system of functions $F'$, such that $F \cup F'$ is a one-ended IFS with the same attractor $A$.*

Finding systems that generate specific attractors is usually quite difficult, and positive results in this direction like [DH92] are few and far between. The question is perhaps easier to approach in a more specific setting:

**Conjecture 5.** *If the attractor $A$ of a piecewise linear IFS $F$ is connected, then there is a system of piecewise linear functions $F'$, such that $F \cup F'$ is a one-ended IFS with the same attractor $A$.*

In the Cayley graph, a walk is determined by a sequence of edges $f \in F$ starting from some vertex. Since that vertex itself is a composition of such functions, we may assume, without loss of generality, that all walks originate from the identity.

**Definition 6.** A walk $(f_1 \circ \ldots \circ f_n)_{n \in \mathbb{N}}$, where each $f_i \in F$, *encodes* a point $x \in M$ if for every point $y \in M$, $\lim_{n \to \infty} f_1 \circ \ldots \circ f_n(y) = x$.

Then the following standard fact about attractors is immediate:

**Lemma 7.** *Every walk encodes a unique point $a \in A$, and conversely for every $a \in A$ there is at least one walk encoding it.*

The next theorem refines the statement of Theorem 3 in that it not only provides a bound on the number of components but establishes a natural correspondence from the ends to the components.

**Theorem 8.** *Suppose the attractor $A$ of an IFS is nondegenerate. Let $\alpha$ denote the map that takes a walk to the point in $A$ it encodes, and let $A'$ be the set of accumulation points in $A$. Then there is a naturally induced map from the set of ends to the components of $A'$, so that the following diagram commutes:*

$$
\begin{array}{ccc}
\text{Rays} & \xrightarrow{\ \alpha\ } & A' \\
\downarrow & & \downarrow \\
\text{Ends} & \dashrightarrow & \mathrm{comp}(A')
\end{array}
$$

(The vertical maps are the respective quotient maps.)

There is a remarkable coincidence of a map in the semigroup being constant, being an idempotent, and being a *dead-end* in the Cayley graph, i.e. a vertex where each originating edge is a loop:

**Lemma 9.** *Suppose $F$ is an IFS. The following conditions for $U \in S(F)$ are equivalent:*

(1) *U is a constant map;*
(2) *There is a map $V \in S$, satisfying $U \circ V = U$;*
(3) *U is an idempotent map, i.e. $U \circ U = U$;*
(4) *U is a dead-end in the Cayley graph $\Gamma(S)$, i.e. $Uf = U$ for all $f \in F$.*

With this in mind, we can explore $\alpha$ on the remainder of its domain.

**Theorem 10.** *Let $C \subseteq A$ denote the image under $\alpha$ of the walks that are not rays. This countable set $C$ consists precisely of all the images of idempotents (i.e. constant maps) in the semigroup. $C$ contains all isolated points of the attractor.*

**Diagram 11.** Indeed, the information from the two theorems can be combined in a bigger commutative diagram that juxtaposes the correspondence between the attractor and its components with that between walks and ends. Here comp denotes connected components; with a subscript it denotes components of a subspace in a larger space.

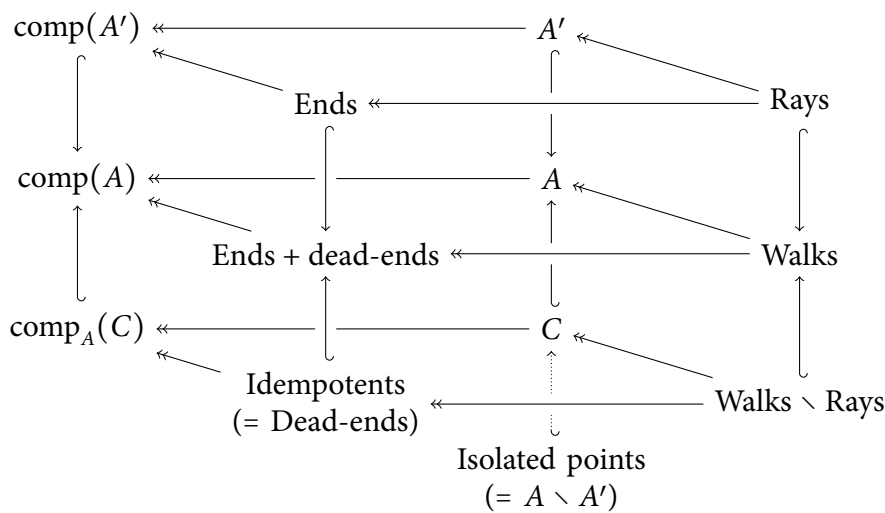

It is linguistic serendipity that the two unrelated notions of ends and dead-ends act so nicely in concert to describe the components. Notice that all vertical maps are simple set inclusions and any object on the middle level is always the union of the ones above and below it, the two in the "front" being a disjoint union. The diagonal maps on the right are different restrictions of $\alpha$.

We now offer a criterion that makes it easy to confirm that various examples given later in the paper are one-ended.

**Definition 12.** Two generators $f, g$ in an IFS $F$ are called *linked*, if there exist $u, v \in S(F)$ such that $f \circ u = g \circ v$. The graph with vertices $F$ and edges $(f, g)$ for all linked pairs is called the *link graph*.

**Theorem 13.** *If $S(F)$ contains no idempotent elements, the link graph of $F$ is connected if and only if $F$ has one end.*

**Example 14.** *In the presence of idempotent elements, neither implication holds.*

In $\mathbb{R}^2$ consider the system of two functions

$$a : (x, y) \mapsto (\tfrac{x}{2}, 0), \quad b : (x, y) \mapsto (0, \tfrac{y}{2}).$$

Note that the attractor is simply $\{0\}$, the fix point of both functions. Evidently, any word in $F$ that is not a power of $a$ or $b$ evaluates to the constant map to 0. In particular the relation $ab \equiv ba$ holds, so the system is linked. It has however two ends.

On the other hand, let

$$a : x \mapsto \tfrac{x}{2}, \quad b : x \mapsto 1$$

be two real functions constituting an IFS, whose attractor is $\{2^{-n} : n \in \mathbb{N}\} \cup \{0\}$. This satisfies exactly the relations $a^k bu \equiv a^k bv$ for $k \in \mathbb{N}$ and any $u, v$ in $S$. It has one end, but is not linked. Note however, that the dead-ends in the Cayley graph coincide with the isolated points of the attractor, as predicted by Theorem 3.

## 3. Further Examples

Most systems defining classical fractals have not just one end, their Cayley graph rather looks like a tree. In this section we will demonstrate how to extend some systems with additional maps to make them one-ended, lending support to Conjecture 4. In the Cayley graph this corresponds to tying the emanating rays in the tree together by virtue of the new maps.

**Example 15.** *The Koch curve, the Sierpiński triangle, the Sierpiński carpet, and the Menger sponge are all attractors of one-ended function systems.*

Consider the Koch curve generated by the two affine functions

$$a : (x, y) \mapsto \left( -\tfrac{x}{2} + \tfrac{y}{2\sqrt{3}} + \tfrac{1}{2}, \ -\tfrac{x}{2\sqrt{3}} - \tfrac{y}{2} + \tfrac{1}{2\sqrt{3}} \right),$$
$$b : (x, y) \mapsto \left( \tfrac{x}{2} + \tfrac{y}{2\sqrt{3}} + 1, \ -\tfrac{x}{2\sqrt{3}} + \tfrac{y}{2} \right).$$

Without changing the attractor, we can add the third function

$$c : (x, y) \mapsto \tfrac{1}{3}(x, y) + \left( \tfrac{1}{3}, \tfrac{2}{9} \right),$$

and simple calculations show that these maps satisfy the relations

$$aab = ca\,,$$
$$cb = bba.$$

By Theorem 13, the system $\{a, b, c\}$ is one-ended.

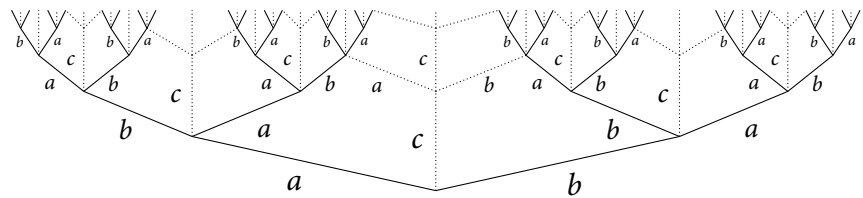

FIGURE 1. The augmented Cayley graph of the Koch curve. The solidly drawn binary tree corresponds to the infinitely-ended free monoid generated by $a$ and $b$ only.

Next consider the Sierpiński triangle, generated by the three functions

$$a : (x, y) \mapsto \left(\tfrac{x}{2}, \tfrac{y}{2}\right), \quad b : (x, y) \mapsto \left(\tfrac{x}{2} + \tfrac{1}{2}, \tfrac{y}{2}\right), \quad c : (x, y) \mapsto \left(\tfrac{x}{2} + \tfrac{1}{4}, \tfrac{y}{2} + \tfrac{\sqrt{3}}{4}\right).$$

Now add two maps contracting the fractal to a line passing through the intersection of two neighbouring triangles,

$$d : (x, y) \mapsto \left(\tfrac{x}{8} + \tfrac{1}{4}, 0\right), \quad e : (x, y) \mapsto \left(\tfrac{x}{8} + \tfrac{y\sqrt{3}}{8} + \tfrac{3}{4}, -\tfrac{x\sqrt{3}}{8} + \tfrac{3y}{8} + \tfrac{\sqrt{3}}{4}\right).$$

These maps satisfy the relations: $abd = dad$, $dbd = bad$, $bce = ebe$, and $ece = cbe$. This gives rise to a spanning tree in the link graph, so by Theorem 13 the system $\{a, b, c, d, e\}$ is one-ended.

Finally, the Sierpiński carpet can be generated by the eight functions $a_1, \ldots, a_4, a_6, \ldots, a_9$ (following the dials on a phone), where each is simply scaling the figure to a third of the size. Additionally, we have two functions $w, e$, corresponding to orthogonal projections to a vertical line and scaling,

$$w : (x, y) \mapsto \left(0, \tfrac{2y}{3} + \tfrac{1}{6}\right), \quad e : (x, y) \mapsto \left(1, \tfrac{2y}{3} + \tfrac{1}{6}\right).$$

Then the following relations hold: $a_1 e = a_2 w$, $a_2 e = a_3 w$, $a_7 e = a_8 w$, $a_8 e = a_9 w$, and also $a_4 w = w a_4$, $a_6 e = e a_6$. Further, $w a_1 a_4 = a_1 a_7 w$, $w a_7 a_4 = a_7 a_1 w$, and $e a_1 a_4 = a_3 a_9 e$. These again satisfy the conditions of Theorem 13, and the system $\{a_1, \ldots, a_4, a_6, \ldots, a_9, w, e\}$ is one-ended.

The same method can be used to construct a one-ended system for the Menger sponge by adding four maps corresponding to four orthogonal projections to parallel edges of the surrounding cube.

Note that all of the above systems consist of affine functions. The next example shows, this is not always possible:

**Example 16.** Consider a crooked version of the Koch curve in the complex plane $\mathbb{C}$, given by the maps $a : z \mapsto \tfrac{z}{3}$, $b : z \mapsto \tfrac{1+i}{3}z + \tfrac{1}{3}$, $c : z \mapsto -\tfrac{i}{3}z + \tfrac{2+i}{3}$, and $d : z \mapsto \tfrac{z}{3} + \tfrac{2}{3}$.

The convex hull of the attractor $A$ is bounded by the lines $K : y = 0$, $L : y = \tfrac{3}{4}x$, $M : y = \tfrac{1}{3}$, and $N : y = -x + 1$. Notice that $K$ intersects $A$ in a Cantor set, while $L$, $M$, and $N$ intersect $A$ in Fort-spaces (countable sets containing exactly one accumulation point), respectively. Further, the accumulation points of two such neighbouring lines never coincide.

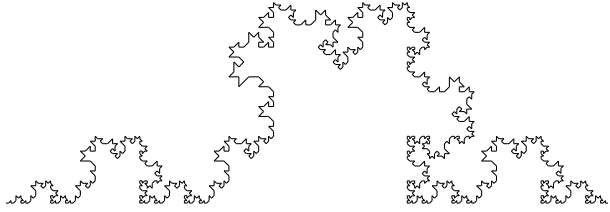

FIGURE 2. A crooked version of the familiar Koch curve that does not allow a one-ended affine system.

For a given planar set $A$, we call a subset of a line an extremal set, if the spanned line segment is part of the boundary of the convex hull of $A$. Then the next geometric fact follows easily:

**Lemma 17.** *An affine transformation takes an extremal set to an extremal set.*

With this we can prove the next statement:

**Theorem 18.** *The affine nontrivial contractions of the fractal $A$ from Example 16 are exactly those generated by words in $a, b, c, d$. Thus there is no affine one-ended IFS generating $A$.*

*Proof.* Let $\varphi$ be an affine map from $A$ to itself. Since $A$ is a curve, $\varphi(A)$ is a segment of $A$ exactly determined by its endpoints $\{p, q\} = \{\varphi(0), q = \varphi(1)\}$. We may assume that $p$ is to the left of $q$ (with respect to the parametrization of the curve), and not both are in the same set $a(A), b(A), c(A),$ or $d(A)$.

Case 1: $p \in a(A) \smallsetminus b(A)$, $q \in b(A) \cup c(A) \cup d(A) \smallsetminus a(A)$. Then $\varphi(A)$ contains a Cantor set in the interval $\left[\frac{1}{3} - \varepsilon, \frac{1}{3}\right]$ for some small $\varepsilon$ that is extremal. Thus $p, q \in \varphi(A \cap K) \subseteq A \cap K$. Since $\varphi$ restricted to $A \cap K$ is affine, this is only possible for $p = 0$ and $q = 1$. Then $\varphi$ has to be the identity and thus is not a contraction.

Case 2: $q \in d(A) \smallsetminus c(A)$, $p \in a(A) \cup b(A) \cup c(A) \smallsetminus d(A)$, works symmetrically to case 1 with a Cantor set situated by the point $\frac{2}{3}$.

Case 3: $p \in b(A) \smallsetminus c(A)$, $q \in c(A) \smallsetminus b(A)$. Then $\varphi(A)$ contains extremal Fort-spaces in $b(A \cap N)$ and $c(A \cap L)$, both with accumulation point $\frac{2+i}{3}$. A contradiction.

This confirms that there are no affine nontrivial contractions of $A$ but those given by words in the generators. Since the semigroup $S = \langle a, b, c, d \rangle$ is free, it cannot be one-ended, and neither can any of its subsemigroups that generate $A$.                                                                    □

However, there is a non-affine one-ended IFS generating the crooked Koch curve. The piecewise linear map $e$ that takes $a(A)$, $b(A)$, $c(A)$, $d(A)$ to $ac(A)$, $ad(A)$, $ba(A)$, and $bb(A)$, respectively, establishes a link between $a$ and $b$, for these now satisfy the relations $ac = ea$ and $ed = bb$. Repeat this method for $b$, $c$ and $c$, $d$; the thus augmented system has one end.

**Example 19.** Generally, the Cayley graph does not necessarily tell the whole story about the attractor. On the real line define the four maps

$$a : x \mapsto \tfrac{x}{2}, \quad b : x \mapsto \tfrac{x}{2} + \tfrac{1}{2}, \quad c : x \mapsto 1, \quad d : x \mapsto \sqrt{2}.$$

The system $\{a, b\}$ has the unit interval as its attractor, as does $\{a, b, c\}$, with the constant map not adding any information. The system $\{a, b, d\}$ shares an isomorphic Cayley graph with the latter (Figure 3), but adds the isolated points $\frac{1+\sqrt{2}}{2}, \frac{3+\sqrt{2}}{4}, \frac{7+\sqrt{2}}{8}, \ldots$ to the attractor.

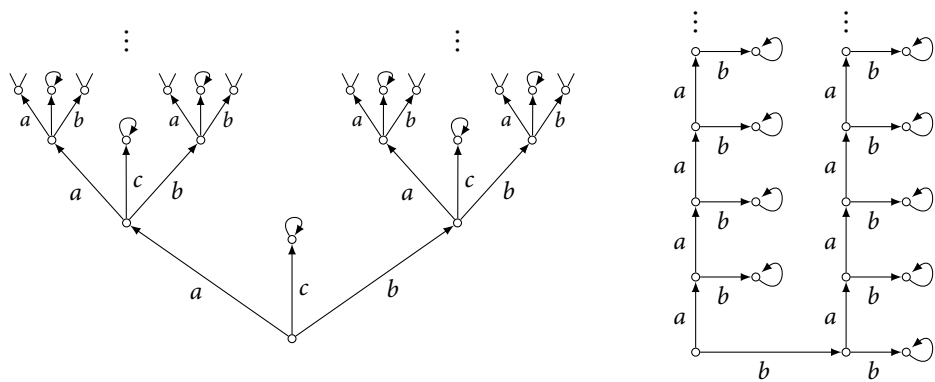

FIGURE 3. (left) The Cayley graph of the two systems generating different attractors from Example 19; (right) the graph from Example 21.

The last two examples shed some light on the various intricate ways idempotents might unexpectedly appear in the semigroup without being evident directly from the iterated function system.

**Example 20.** In $\mathbb{R}^k$ of arbitrary dimension $k$ consider the piecewise defined map

$$a : x \mapsto \begin{cases} 0 & \text{if} \quad |x| < 2, \\ x/2 - x/|x| & \text{if} \quad 2 \le |x| < 2^n, \\ (2^{n-1} - 1)x/|x| & \text{if} \quad 2^n \le |x|, \end{cases}$$

for some $n > 1$. Obviously the attractor of the IFS made up of this single map has to be just the origin. The Cayley graph however has precisely the $n$ distinct vertices $a^1, \ldots, a^n$. The last, $a^n$, is a constant function and thus a dead-end with only a loop going out from it.

**Example 21.** As a variation of Example 14, consider the system of two functions

$$a : (x, y) \mapsto \left(\tfrac{x}{2}, 0\right), \quad b : (x, y) \mapsto \left(1, \tfrac{x}{2}\right)$$

in the plane. Neither is itself constant. The attractor is a countable set with two accumulation points, precisely, $A = \{(0,0)\} \cup \{(2^{-n}, 0) : n \geq 0\} \cup \{(1, 2^{-n}) : n > 0\}$. The Cayley graph, with the two ends and the infinitely many dead-ends clearly visible, is depicted in Figure 3.

## 4. Proofs

**Lemma 22.** *Let $a, b \in X$ be two points in a compact metric space. If for every $\varepsilon > 0$ there is a sequence of points $a = c_1, \ldots, c_n = b \in X$ such that the distance $d(c_i, c_{i+1})$ between consecutive elements is less than $\varepsilon$, then $a$ and $b$ are in the same connected component of $X$.*

*Proof.* Suppose $a$ and $b$ are in distinct components of $X$. In a compact Hausdorff space, the components coincide with the quasi-components, thus $X$ is the disjoint union of open sets $A$ and $B$ containing $a$ and $b$, respectively. Then $A$ and $B$ are also compact and $\varepsilon := d(A, B) > 0$. Pick a sequence $a = c_1, \ldots, c_n = b \in X$ as in the hypothesis. For some index $i$, $c_i \in A$ and $c_{i+1} \in B$. But then $d(c_i, c_{i+1}) \geq \varepsilon$, a contradiction. $\square$

*Proof of Lemma 9.* Suppose $U \in S$ is the constant map to $x \in M$. Then so is $U \circ U$, and they are equal. If on the other hand, $U \circ V = U$, then by recursion $U = U \circ V^n$. So for arbitrary points $x, y \in M$ and $\lambda$ the contraction factor of $F$ (i.e. $d(x, y) \leq \lambda d(f(x), f(y))$ for all $f \in F$ and all $x, y \in M$), the distance of their images $d(U(x), U(y)) = d(U \circ V^n(x), U \circ V^n(y)) \leq \lambda^{n+1} d(x, y)$ has to be zero, thus $U$ has to be a constant map. $\square$

Next we are going to interpret different representations of a single element in $S$ as elements in the free semigroup with alphabet $F$. For that purpose, call any finite sequence of maps in $F$ *a word*, written as $f_1 \ldots f_n$. The *length of a word* is then the length of the sequence. Two words are *congruent*, $f_1 \ldots f_n \equiv g_1 \ldots g_m$, if the corresponding evaluations in $S$ agree: $f_1 \circ \ldots \circ f_n = g_1 \circ \ldots \circ g_m$. Lastly, let $B_k$ denote the ball of radius $k$ in $S$, i.e. the set of maps that can be represented by a word of length at most $k$.

**Lemma 23.** *If a walk $(f_1 \circ \ldots \circ f_n)_{n \in \mathbb{N}}$ is not a ray in $\Gamma(S)$, then, for some $k$, $f_1 \circ \ldots \circ f_k$ is equal to the constant map to the point $a \in A$ encoded by the sequence.*

*On the other hand, if there is no constant map in $S$, then for each $k \in \mathbb{N}$ there is $N \in \mathbb{N}$ so that no word of length at least $N$ can be replaced by a word of length at most $k$, i.e.*

$$\forall k \in \mathbb{N} \; \exists N \in \mathbb{N} \; \forall n \geq N \; \forall (f_i)_{i=1}^n : f_1 \circ \ldots \circ f_n \in S \smallsetminus B_k.$$

*Proof.* If the sequence is not a ray, two terms in it are not distinct. But then $f_1 \circ \ldots \circ f_k = f_1 \circ \ldots \circ f_k \circ f_{k+1} \circ \ldots \circ f_n$ for some $n, k \in \mathbb{N}$, so by Lemma 9, $f_1 \circ \ldots \circ f_k$ is constant.

Next assume the second claim were false. By the pigeonhole principle, since there are only finitely many elements in $S$ that can be represented by words of length $k$, there is $U \in S$ that can be represented by words of arbitrarily large length. But then $d(U(x), U(y)) \leq \lambda^n d(x, y)$ for points $x, y \in M$ and arbitrary $n \in \mathbb{N}$, so $U$ would be a constant map, a contradiction. $\square$

**Lemma 24.** *Each $a \in A'$ is encoded by some ray.*

*Proof.* Let $a \in A'$. Then there exists a sequence $(a_k)_{k \in \mathbb{N}}$ of points in $A$ converging to $a$. For each $a_k$ let $(f_{k,1} \circ \ldots \circ f_{k,n})_{n \in \mathbb{N}}$ be an encoding sequence. There is an infinite subsequence $(j_1(k))_{k \in \mathbb{N}}$ of $(k)_{k \in \mathbb{N}}$ such that all $f_{j_1(k),1}$ are equal independently of $k$. Proceeding inductively, there is a subsequence $(j_n(k))_{k \in \mathbb{N}}$ of $(j_{n-1}(k))_{k \in \mathbb{N}}$ such that all $f_{j_n(k),n}$ are equal. Then the diagonalized sequence $(f_{j_1(1),1} \circ \ldots \circ f_{j_n(1),n})_{n \in \mathbb{N}}$ has to be a ray since the image of each term contains $\{a_{j_n(k)} : k \in \mathbb{N}\}$ and is therefore infinite and not constant. Further, as the images are closed and their intersection consists of exactly $a$, the sequence encodes $a$. $\square$

**Lemma 25.** *If two rays are equivalent, the encoded points are in the same connected component of the attractor.*

*Proof.* Let $a, b$ be points encoded by rays $(f_1 \circ \ldots \circ f_n)_{n \in \mathbb{N}}$, $(g_1 \circ \ldots \circ g_n)_{n \in \mathbb{N}}$, respectively. Let $D$ be the diameter of the attractor and $\lambda$ be the contraction factor of the IFS; then for $\varepsilon > 0$ pick $k \in \mathbb{N}$ such that $\lambda^k D < \varepsilon$ and let $B_k$ be the set of all elements in $S$ that can be represented by words in $F$ of length less than $k$. Because the above are rays, there is an index $\ell \in \mathbb{N}$, such that $f_1 \circ \ldots \circ f_\ell$ and $g_1 \circ \ldots \circ g_\ell$ are not in $B_k$. Since the sequences are equivalent and $B_k$ is finite, there is a path $s_1, \ldots, s_n \in S \smallsetminus B_k$ such that

$$s_1 = f_1 \circ \ldots \circ f_\ell, \quad s_n = g_1 \circ \ldots \circ g_\ell,$$
$$\text{for each } i \text{ there is } f \in F \; : \; s_{i+1} = s_i \circ f \text{ or } s_i = s_{i+1} \circ f.$$

Choose $f_0 \in F$ nonconstant and $y \in A$ arbitrarily. Then the points $c_i := \lim_{n \to \infty} s_i \circ f_0^n(y)$ are encoded by rays and satisfy $d(c_i, c_{i+1}) \leq \lambda^k D < \varepsilon$, that is the condition of Lemma 22. Hence $a$ and $b$ are in the same connected component. $\square$

*Proof of Theorem 8.* First suppose the sequence $(f_1 \circ \ldots \circ f_n)_{n \in \mathbb{N}}$ encodes an isolated point $a \in A \smallsetminus A'$. There is an $\varepsilon$-neighbourhood containing no other point of the attractor. Let $D$ be the diameter of $A$ and $\lambda$ the contraction factor of $F$. Then for some $n$, the diameter of the image of the attractor under $f_1 \circ \ldots \circ f_n$ is bound by $\lambda^n D < \varepsilon$, so either the whole attractor consists of just a single point or that map is constant to $a$. But according to Lemma 9 a constant map corresponds to a dead-end in the graph, and thus the walk $(f_1 \circ \ldots \circ f_n)_{n \in \mathbb{N}}$ cannot be a ray. Thus $A'$ includes $\alpha(\text{Rays})$, and due to Lemma 24, $\alpha(\text{Rays})$ conversely includes $A'$. So the map on top of the diagram is well-defined and surjective. By Lemma 25 two equivalent rays encode points which lie in the same component of $A$ (and thus also of $A'$), so the induced map is well-defined. Consequently, the diagram commutes.  □

*Proof of Theorem 10.* Given a walk $w$ that is not a ray, Lemma 23 associates to $w$ the constant map to the point $\alpha(w)$. As the semigroup is countable, so is the number of these constant maps and hence so is $C$.

Since $\alpha(\text{Rays})$ is included in $A'$, the isolated points must be encoded by walks that are not rays and thus are included in $C$.  □

*Proof of Theorem 3.* A nondegenerate component has to be a component in $A'$. Then the first claim follows immediately from the fact that the map in Theorem 8 is surjective. The isolated points are contained in the set $C$ from Theorem 10, and that countable set corresponds to the idempotent elements by Lemma 9.  □

*Proof of Theorem 13.* To show $F$ is one-ended amounts to finding for each $k \in \mathbb{N}$ some $N \in \mathbb{N}$, so that each pair of elements in $S$ outside the ball of radius $N$ can be connected by a path in the Cayley graph $\Gamma(S)$ that lies outside the ball of radius $k$. Lemma 23 guarantees that the evaluation of any word of length $N$ is in $S \smallsetminus B_k$. An induction argument will show that words of length $n$ can be connected by a sequence of words of length at least $n$. Since longer words can be shortened by a path to a word of length $n$, this will prove one direction.

Consider the case $n = 1$. Let $f, g$ be in $F$. Since the link graph is connected, there is a sequence $f = w_1, w_2, \ldots, w_k = g$ in $F$ so that for consecutive elements there are words $u_i, v_i$ so that $w_i u_i \equiv w_{i+1} v_i$ (for $1 \leq i < k$). Clearly, there is a path from $w_i$ to $w_i u_i$ and then from $w_{i+1} v_i$ to $w_{i+1}$, each of them consisting of words of length at least 1. Splicing all these together gives the required path from $f$ to $g$.

Now let $fU$ and $gV$ be words of length $n > 1$. Again there is a sequence $f = w_1, w_2, \ldots, w_k = g$ in $F$ so that for consecutive elements there are words $u_i, v_i$ so that $w_i u_i \equiv w_{i+1} v_i$. By extending the $u_i$ and $v_i$, if necessary, with arbitrary words, we may assume all of them have length at least $n - 1$. By induction there is a path from $U$ to $u_1$, for each $i$ one from $v_i$ to $u_{i+1}$, and one from $v_k$ to $V$, all of them consisting of words of length at least $n - 1$. Multiplication with $f$, $w_{i+1}$, and $g$ from the left then yields paths from $fU$ to $fu_1$, from $w_{i+1}v_i$ to $w_{i+1}u_{i+1}$ and from $gv_k$ to $gV$, respectively. Again splicing all these together results in a path from $fU$ to $gV$ consisting of words of length at least $n$.

Concerning the other direction, first notice that, since there are no idempotent elements in $S$, the orbit $\{f^k : k \in \mathbb{N}\}$ of a single generator is always an infinite set of vertices in the Cayley graph. Assuming $F$ has only one end then implies that for two $f, g \in F$ there is some exponent $n \in \mathbb{N}$, so that both $f^n$ and $g^n$ are connected in $S \smallsetminus B_0$. Consider a path between these, and further a sequence of words in $F$ evaluating to its vertices. We may assume the first word to be $f^n$, the last to be $g^n$. Two consecutive words in the sequence correspond to endpoints of an edge and thus evaluate to $u$ and $u \circ h$, or vice versa, for some $u \in S$, $h \in F$. Hence we have a congruence of two words representing $u \circ h$, which shows the first symbols of the two words are linked. Following along the sequence of words then also determines a path in the link graph between $f$ and $g$, and since these were chosen arbitrarily, the link graph is connected. $\qquad\square$